\documentclass[11pt]{amsart}
\usepackage{a4wide}
\usepackage{amsmath, amsfonts, amssymb,amscd,graphicx,amsthm}
\usepackage{mathrsfs,url}

\date{\today}

\author
{Manuel Bodirsky}
\address{Institut f\"{u}r Algebra\\TU Dresden\\01062 Dresden\\Germany}
    \email{Manuel.Bodirsky@tu-dresden.de}
   \urladdr{http://www.math.tu-dresden.de/~bodirsky/}
    \thanks{The first author has received funding from the European Research Council (Grant Agreement no.~681988, CSP-Infinity). 
The second author has received funding from the  Austrian Science Fund (FWF) through projects I836-N23 and P27600,
and from the Czech Science Foundation (grant No 18-20123S).
The first and third author have received funding from the European Research Council under the European Community's Seventh Framework Programme (FP7/2007-2013 Grant Agreement no. 257039). 
The third author has been supported by the Hungarian Science Fund (OTKA) grant no.~K109185, by the EFOP-3.6.2-16-2017-00015 project, which has been supported by the European Union, co-financed by the European Social Fund, by the National Research, Development and Innovation Fund of Hungary, financed under the FK 124814 and PD 125160 funding schemes, the J{\'a}nos Bolyai Research Scholarship of the Hungarian Academy of Sciences, and by the \'{U}NKP-18-4 New National Excellence Program of the Ministry of Human Capacities.
}

\author
{Michael Pinsker}
	\address{Institut f\"{u}r Diskrete Mathematik und Geometrie, FG Algebra, TU Wien, Austria, and Department of Algebra, Charles University, Czech Republic}        
	\email{marula@gmx.at}
    \urladdr{http://dmg.tuwien.ac.at/pinsker/}

\author
{Andr\'{a}s Pongr\'{a}cz}
    \address{Department of Algebra and Number Theory\\
    University of Debrecen\\
    4032 Debrecen, Egyetem square 1\\
    Hungary}
\email{pongracz.andras@science.unideb.hu}

\title{Projective clone homomorphisms}

\newcommand{\fC}{\mathscr C}

\DeclareMathOperator{\typ}{typ}

\newcommand{\ignore}[1]{}

\newcommand{\cl}[1]{\langle #1 \rangle}

\newcommand{\To}{\rightarrow}

\DeclareMathOperator{\Csp}{CSP}

\DeclareMathOperator{\Aut}{Aut}

\DeclareMathOperator{\Pol}{Pol}

\newcommand{\F}{\mathscr F}

\theoremstyle{plain}

    \newtheorem{thm}{Theorem}[section]

    \newtheorem{lem}[thm]{Lemma}

    \newtheorem{prop}[thm]{Proposition}

    \newtheorem{quest}[thm]{Question}
    
    \newtheorem{conj}[thm]{Conjecture}

\theoremstyle{definition}

    \newtheorem{defn}[thm]{Definition}

\DeclareMathOperator{\Pro}{\mathscr P}

\newcommand{\red}[1]{#1}

\begin{document}
\begin{abstract}
It is known that a countable $\omega$-categorical structure interprets all finite structures primitively positively if and only if its polymorphism clone maps to the clone of projections on a two-element set via a continuous clone homomorphism. We investigate the relationship between the existence of a clone homomorphism to the projection clone,  and the existence of such a homomorphism which is continuous and thus meets the above criterion.
\end{abstract}

\maketitle

\section{Introduction}
\label{sect:intro}
A \emph{function clone} is a set of finitary functions on a fixed domain which is closed under composition and which contains all projections. There are two main sources of function clones: the set of \emph{term operations} 
of any algebra $\mathfrak A$ is a function clone, and in fact every function clone is of this form; moreover, the \emph{polymorphism clone} $\Pol(\Gamma)$ of any first-order structure $\Gamma$, consisting of all finitary functions which \emph{preserve} $\Gamma$, forms a function clone. Reminiscent of automorphism groups and endomorphism monoids, polymorphism clones carry information about the structure that induces them, and are a powerful tool in the study of first-order structures. And  similarly to the situation for permutation groups and transformation monoids, a function clone is a polymorphism clone if and only if it is \emph{closed} in the topology of \emph{pointwise convergence} on the set of all finitary functions on a fixed domain. This topology is obtained by viewing the domain as a discrete space, and equipping for all $n\geq 1$ the $n$-ary functions on it with the product topology; finally, the set of all finitary functions is the sum space of these spaces (cf.~for example~\cite{Topo-Birk, BP-reductsRamsey}). On countable domains, this topology is induced by the following metric. 
The distance of two functions of different arity is 1, and that of distinct functions of the same arity $n$ is $2^{-k}$, where $k$ is the smallest index where the two functions differ, in a fixed enumeration of the set of $n$-tuples of elements of the domain.
In the sequel, we shall simply say that a function clone is closed iff it is closed in this topological space, and Cauchy sequences are meant with respect to this metric space. 

Again bearing analogy to automorphism groups, abstract properties of the polymorphism clone of a structure or of the term clone of an algebra can translate into properties of the structure or algebra, respectively. 
Such abstract properties can be purely \emph{algebraic}, e.g., given by equations which hold in the function clone, or  \emph{algebraic and topological}, i.e., 
captured if we consider in addition the topological structure of a clone given by pointwise convergence. In the analogy with permutation groups where the latter kind of abstraction leads to the notion of a topological group, we here obtain \emph{topological clones}~\cite{Reconstruction}.

On every set there is a smallest function clone, namely the function clone which consists precisely of the finitary projections over this set. For any two sets  of at least two elements, these projection clones are isomorphic algebraically and topologically; in fact, the topology on any such projection clone is discrete. We write $\Pro$ for the topological clone induced by any such projection clone, and denote for $1\leq k\leq n$ its $n$-ary projection to the $k$-th coordinate by $\pi^n_k$. 
Given a function clone $\fC$, it is an important structural property whether or not there exists a \emph{continuous clone homomorphism} to the clone $\Pro$, i.e., a continuous  mapping $\xi\colon \fC\To\Pro$ which 
\begin{itemize}
\item sends functions of $\fC$ to functions of the same arity in $\Pro$,
\item sends projections in $\fC$ to the corresponding projections in $\Pro$, and
\item preserves composition, i.e., for all $f,g_1,\ldots,g_n\in\fC$
$$\xi(f(g_1,\ldots,g_n))=\xi(f)(\xi(g_1),\ldots,\xi(g_n)).$$
\end{itemize}
We call clone homomorphisms to $\Pro$ \emph{projective}.
For countable $\omega$-categorical structures (cf.~\cite{Hodges} for standard model theoretic definitions), the importance of the existence of a continuous projective homomorphism is given by the following observation. Here, a primitive positive interpretation of a structure $\Delta$ in a structure $\Gamma$ is an interpretation in the sense of classical model theory (cf.~\cite{Hodges}) in which all involved formulas are \emph{primitive positive}, i.e., existentially quantified conjunctions of atomic formulas.

\begin{prop}[Bodirsky and Pinsker~\cite{Topo-Birk}]\label{prop:allstructures}
Let $\Gamma$ be a countable $\omega$-categorical structure. Then the following are equivalent:
\begin{itemize}
\item  $\Pol(\Gamma)$ has a continuous projective homomorphism;
\item all finite structures have a primitive positive interpretation in $\Gamma$;
\item the structure $(\{0,1\};\{(1,0,0),(0,1,0),(0,0,1)\})$ has a primitive positive interpretation in $\Gamma$.
\end{itemize}
\end{prop}

When all finite structures have  a primitive positive interpretation in a structure $\Gamma$, then $\Gamma$ can be considered at least as complicated as all finite structures in the quasiorder of primitive positive interpretations on first-order structures. This property is of particular interest in \emph{constraint satisfaction} in theoretical computer science. For any structure $\Gamma$ in a  finite relational language, the \emph{constraint satisfaction problem of $\Gamma$}, denoted by $\Csp(\Gamma)$, is the problem of deciding whether or not a given conjunction of atomic formulas in the language of $\Gamma$ has a solution in $\Gamma$. It is not hard to see from the definitions that when a structure $\Delta$ has a primitive positive interpretation in a structure $\Gamma$, then $\Csp(\Delta)$ is polynomial-time reducible to $\Csp(\Gamma)$. Therefore, if the conditions of Proposition~\ref{prop:allstructures} hold for a structure $\Gamma$, then $\Csp(\Gamma)$ is, up to polynomial time, at least as hard as the constraint satisfaction problem of any finite structure, and in particular NP-hard. 

In many situations, the existence of a continuous projective homomorphism is even believed to be the only possible source of NP-hardness: for example, for finite $\Gamma$, the famous \emph{tractability conjecture} stated that under certain conditions on $\Gamma$ which can be assumed without loss of generality, 
$\Csp(\Gamma)$ is NP-complete if $\Pol(\Gamma)$ has such a homomorphism, and in P otherwise. 
\red{Two different proofs of this conjecture have been announced in~\cite{BulatovFVConjecture} and in~\cite{ZhukFVConjecture}. }
Note that in the finite case, the existence of a continuous projective homomorphism is a purely algebraic property of $\Pol(\Gamma)$, since the topology of any function clone on a finite set is discrete. For a large and natural class of infinite structures $\Gamma$, a similar conjecture has been formulated. A Fra\"{i}ss\'{e} class of finite structures in a finite relational language $\tau$ is called \emph{finitely bounded} iff there exist finitely many finite ``forbidden'' $\tau$-structures such that the class consists precisely of those finite $\tau$-structures which do not embed any of the forbidden structures. Of particular interest in constraint satisfaction are structures which are first-order definable in Fra\"{i}ss\'{e}-limits of finitely bounded Fra\"{i}ss\'{e} classes; such structures are $\omega$-categorical. For every $\omega$-categorical structure $\Gamma$ there is an $\omega$-categorical structure $\Delta$ whose automorphisms are dense in its endomorphisms (in the above topology) and which is \emph{homomorphically equivalent} to $\Gamma$, i.e., $\Gamma$ can be homomorphically mapped into $\Delta$ and vice versa~\cite{Cores-journal, BodHilsMartin,BKOPP, BKOPP-conf}. The structure $\Delta$ is called the \emph{model-complete core} of $\Gamma$, and $\Csp(\Delta)$ has the same true instances as $\Csp(\Gamma)$.

\begin{conj}[Bodirsky and Pinsker, 2012]\label{conj:tractability}
Let $\Gamma$ be the model-complete core of a structure which is first-order definable in the Fra\"{i}ss\'{e}-limit of a  finitely bounded Fra\"{i}ss\'{e} class. Then precisely one of the following holds:
\begin{itemize}
\item there exists an expansion $\Gamma' $of $\Gamma$ by finitely many constants such that $\Pol(\Gamma')$ has a continuous projective homomorphism (and $\Csp(\Gamma)$ is NP-complete by Proposition~\ref{prop:allstructures} and the fact that such expansions do not increase the complexity of the CSP);
\item for any such expansion $\Gamma'$ the clone $\Pol(\Gamma')$ has no continuous projective homomorphism, and $\Csp(\Gamma)$ is in P.
\end{itemize}
\end{conj}

\red{Recent progress on this conjecture has been made by finding equivalent formulations~\cite{wonderland,BKOPP,BKOPP-conf,Topo} and by confirming it in special cases~\cite{BMPP16,MMSNP,BodMot-Unary}.} 
One approach to proving the conjecture would be showing that if there is no continuous projective homomorphism of $\Pol(\Gamma')$, then there is no projective homomorphism at all, and furthermore showing that in that situation $\Csp(\Gamma)$ is in P. 
In this paper, we investigate the first, complexity-free, part. In other words, we investigate the following question:

\begin{quest}\label{quest:main}
Let $\fC$ be a closed function clone with a projective homomorphism. Does $\fC$ also have a continuous projective homomorphism?
\end{quest}

It is worth noting that since the topological clone $\Pro$ is discrete, a mapping to $\Pro$ is continuous if and only if the preimage of each projection $\pi^n_k$ in $\Pro$ is a clopen set. Therefore, a continuous clone homomorphism $\xi\colon \fC\To \Pro$ gives us for every $n\geq 1$ a partition of the $n$-ary functions of $\fC$ into $n$ clopen sets such that composition of representatives of those sets behaves like composition of projections.

We mention a related, and for polymorphism clones of countable $\omega$-categorical structures equivalent, formulation of Question~\ref{quest:main}. 
For an algebra $\mathfrak A$, the \emph{variety generated by $\mathfrak A$} is the class of all algebras that can be obtained from $\mathfrak A$ by taking arbitrary powers, subalgebras, and homomorphic images. The \emph{pseudovariety generated by $\mathfrak A$} is defined similarly, with the only difference being that only \emph{finite} powers are allowed. Given a function clone $\fC$, we can view its functions as the functions of an algebra by giving it a signature in an arbitrary way. In particular, this way we can make an algebra out of any polymorphism clone. 
The notions of a variety and pseudovariety then relate to clone homomorphisms; in particular, we have the following. An algebra is \emph{trivial} iff all of its operations are projections.

\begin{thm}[Birkhoff~\cite{Bir-On-the-structure}; Bodirsky and Pinsker~\cite{Topo-Birk}]
Let $\mathfrak A$ be an algebra whose operations constitute the polymorphism clone $\fC$ of an $\omega$-categorical structure on a countable domain. Then:
\begin{itemize}
\item The variety generated by $\mathfrak A$ contains a two-element trivial algebra if and only if $\fC$ has a projective homomorphism.
\item The pseudovariety generated by $\mathfrak A$ contains a two-element trivial algebra if and only if $\fC$ has a continuous projective homomorphism.
\end{itemize}
\end{thm}

By investigating Question~\ref{quest:main}, we therefore investigate in this paper when certain algebras have trivial algebras in the variety or pseudovariety they generate.

\section{Results}
\label{sect:results}
As a first observation, we note that in order to prove that a homomorphism
$\xi \colon \fC \to \Pro$ is continuous, it suffices to focus on the \emph{binary} functions in the clone, that is, it suffices to show that the preimages 
of the two binary projections under $\xi$ are open (Section~\ref{sect:binary}). 
This is false if we replace $\Pro$ by other clones, even function clones on a finite domain. We moreover show that ``partial homomorphisms'' from binary fragments of function clones to $\Pro$ can always be extended to the function clone they generate, which is again a special property of $\Pro$.

A potential strategy for obtaining a positive answer 
to Question~\ref{quest:main} is to prove something stronger than the property asked in Question~\ref{quest:main}.
Let $\fC$ be a closed function clone that has a homomorphism $\xi$ to $\Pro$. 
Instead of proving that there also exists a continuous homomorphism from $\fC$ to $\Pro$, we might ask whether $\xi$ itself is necessarily continuous. 
\begin{quest}\label{prob:allcont}
Let $\fC$ be a \red{closed} function clone. Is every projective homomorphism of $\fC$ continuous?
\end{quest}
In Section~\ref{sect:auto-cont} we provide two examples of \red{closed} function clones over countable base sets where Question~\ref{prob:allcont} has a negative answer. The two examples can be seen as opposite extreme cases of function clones: one is \emph{oligomorphic}, 
that is, it contains a permutation group which
has for each $n\geq 1$ finitely many orbits in its componentwise action on $n$-tuples. 
The other one is the term clone of a \emph{locally finite} algebra, that is,
the finitely generated subalgebras of the algebra are finite. 
Note that the clone of a locally finite algebra is never oligomorphic. We also remark that the closed oligomorphic function clones on a countable domain are precisely the polymorphism clones of $\omega$-categorical structures on this domain.

However, both examples rely on the existence of non-principal ultrafilters on a countable set. Hence, we cannot exclude the existence of models of ZF where Question~\ref{prob:allcont} has a positive answer. 
It should be mentioned in this context that there are models of ZF+DC (axiom of dependent choice) where every homomorphism between closed subgroups of the full symmetric group on a countable set  is continuous (see the discussion in Section 8 in \cite{Topo-Birk}). Whether similar phenomena prevail in the world of clones
is not known, and we refer to the discussion in~\cite{Reconstruction}.


Question~\ref{quest:main}, on the other hand, might well have a positive answer 
for \emph{all} closed function clones. 
We give a positive answer to Question~\ref{quest:main} 
for all closed clones of idempotent locally finite algebras (Section~\ref{sect:locfin}). Here we use 
recent results about finite idempotent algebras~\cite{Siggers} 
via a compactness argument. 

Again using results from finite idempotent algebras, but this
time in a less obvious way, we also give (in Section~\ref{sect:oligo}) 
a positive answer to Question~\ref{quest:main} for an important class of oligomorphic function clones which we describe next. 
Let $\Delta$ be a structure with a finite relational
signature and domain $D$ which is \emph{homogeneous}, i.e., any isomorphism between finite substructures of $\Delta$ extends to an automorphism of $\Delta$. We say that a function $f \colon D^n \to D$
is \emph{canonical} with respect to $\Delta$ if
for all $k\geq 1$, all $k$-tuples $t_1,\dots,t_n \in D^k$, and for all $\alpha_1,\ldots,\alpha_n \in \Aut(\Delta)$ there exists a $\beta \in \Aut(\Delta)$ such that $f(\alpha(t_1),\dots,\alpha(t_n)) = \beta(f(t_1,\dots,t_n))$ (where we apply functions to tuples componentwise). We prove that Question~\ref{quest:main} has a positive
answer for every closed polymorphism clone containing $\Aut(\Delta)$, and all of whose operations are canonical
with respect to $\Delta$. 
Clones of canonical functions arise naturally
when the finite substructures of a homogeneous structure 
$\Delta$ have the \emph{Ramsey property}
(for details, see~\cite{BP-reductsRamsey}), and are in this case of crucial importance in the study of CSPs of structures which are definable in $\Delta$~\cite{BP-reductsRamsey, BodPin-Schaefer}.
In fact, it can be shown in this situation that every closed function clone that contains 
$\Aut(\Delta)$ 
can be written as $\overline{\bigcup_{i=1}^{\infty} \fC_i}$ 
where each $\fC_i$ is a function clone consisting of operations that
are canonical with respect to an expansion of $\Delta$ 
with finitely many constants. This fact, and canonical operations in general,
have served as the main tool in the successful verifications of Conjecture~\ref{conj:tractability} so far. 


The condition that $\fC$ be closed cannot be omitted from Question~\ref{quest:main}. We present a non-closed, oligomorphic counterexample 
in Section~\ref{sect:oligo}.

We close with a list of open problems that are related to our two research questions (Section~\ref{sect:open}).


\section{The Binary Fragment}
\label{sect:binary}

\begin{defn}
Let $\F$ be a set of finitary functions on some set. The \emph{function clone generated by $\F$}, denoted by $\cl{\F}$, is the smallest function clone containing $\F$. It consists of all {term functions} over $\F$, i.e., functions that are obtained by composing functions from $\F$ and projections. The  \emph{closed function clone generated by $\F$}, denoted by $\overline{\cl{\F}}$, is the smallest closed function clone containing $\F$.  It is the closure of $\cl{\F}$ in the space of finitary functions on the domain, and consists of those functions which agree with some function in $\cl{\F}$ on every finite set.
\end{defn}

\begin{defn}
Let $\F$ be a set of finitary functions on some set (not necessarily a function clone). A \emph{projective partial homomorphism} of $\F$ is a mapping from $\F$ to $\mathbf 1$ which preserves arities of functions, which sends every projection in $\F$ to the same projection in $\mathbf 1$, and which preserves composition whenever it is defined in $\F$.
\end{defn}

\begin{lem}\label{lem:extendfrombinary}
Let $\F$ be the binary fragment of a function clone, and let $\xi$ be a projective partial homomorphism of $\F$. Then $\xi$ extends uniquely to a projective homomorphism of the function clone $\cl{\F}$ generated by $\F$.
\end{lem}
\begin{proof}
Let $\fC$ be the function clone generated by $\F$. It is clear that any extension of $\xi$ to a projective homomorphism of $\fC$ is unique. To see that such an extension exists, define for any $n\geq 1$, any $n$-ary $f\in\fC$, and any $1\leq i\leq n$ a binary function $f_i(x,y):=f(x,\ldots,x,y,x,\ldots,x)$, where the variable $y$ is inserted at the $i$-th coordinate. Clearly, $f_i\in\F$. By a straightforward induction on terms over $\F$, one can show that there exists precisely one $1\leq i\leq n$ for which $\xi(f_i)$ equals $\pi_2^2$. We set $\xi'(f):=\pi^n_i$ for that particular $i$. It is easy to verify that defined this way, $\xi'$ is a homomorphism from $\fC$ to $\mathbf 1$. Clearly, $\xi'$ extends $\xi$.
\end{proof}

\begin{lem}\label{lem:continuousbinary}
Let $\fC$ be a function clone, and let $\xi\colon \fC\To\mathbf 1$ be a homomorphism. Then $\xi$ is continuous if and only if its restriction to the binary fragment of $\fC$ is.
\end{lem}
\begin{proof}
	For $1\leq i\leq n$, and an $n$-ary $f\in\fC$, we define $f_i$ as in the proof of Lemma~\ref{lem:extendfrombinary}. 
	 Then $\xi^{-1}(\{\pi^n_i\})$ consists of those $n$-ary functions $f$ for which $f_i\in \xi^{-1}(\{\pi^2_2\})$, a clopen set since $\xi^{-1}(\{\pi^2_2\})$ is clopen and since the mapping which sends every $n$-ary $f\in\fC$ to $f_i$ is continuous.
\end{proof}

\begin{lem}\label{lem:extendfrombinarycont}
Let $\F$ be the binary fragment of a closed function clone \red{on a countable domain}, and let $\xi$ be a continuous partial projective homomorphism of $\F$. Then $\xi$ extends uniquely to a continuous projective homomorphism of the closed function clone $\overline{\cl{\F}}$ generated by $\F$.
\end{lem}
\begin{proof}
Let $\fC'$ be the function clone generated by $\F$. By Lemma~\ref{lem:extendfrombinary}, $\xi$ extends uniquely to  a projective homomorphism  $\xi'$ of $\fC'$, and this extension is continuous by Lemma~\ref{lem:continuousbinary}. Since the closed function clone $\fC$ generated by $\F$ is just the completion of the topological clone $\fC'$, we have to show that $\xi$ is Cauchy continuous in order to prove that $\xi'$ extends continuously to $\fC$. Let $(f^j)_{j\in\omega}$ be a Cauchy sequence of $n$-ary functions in $\fC'$, where $n\geq 1$. Then for all $1\leq i\leq n$, the sequences $(f^j_i)_{j\in\omega}$ are Cauchy as well, and hence they converge. By the continuity of $\xi$, their value under $\xi$ is constant for all $j\geq k$, for some $k\in\omega$. Hence, there exists $1\leq i\leq n$ such that $\xi'(f^j)=\pi^n_i$ for all $j\geq k$, showing Cauchy continuity. It is straightforward to check that the continuous extension to $\fC$ is a homomorphism.
\end{proof}

\section{Automatic Continuity to $\Pro$}
\label{sect:auto-cont}

Recall that an $n$-ary operation $f$ on a set $D$ is called \emph{conservative} iff $f(x_1,\ldots,x_n)\in\{x_1,\ldots,x_n\}$ for all $x_1,\ldots,x_n\in D$. Note that any conservative function $f$ on $D$ is  \emph{idempotent}, i.e., $f(x,\ldots,x)=x$ for all $x\in D$.

\begin{prop}\label{prop:discont1}
There exists a function clone of conservative functions on a countable set with a discontinuous projective homomorphism.
\end{prop}
\begin{proof}
Let $\F$ be the set of all binary functions $f$ on $\omega$ such that $\{f(a,b),f(b,a)\}=\{a,b\}$ for all $a,b\in\omega$. Then all functions in $\F$ are conservative, and for all $a,b\in\omega$ the restriction of $f$ to $\{a,b\}$ equals the restriction of a binary projection to this set. Let $\fC:=\cl{\F}$. Then the binary functions of $\fC$ are precisely the functions in $\F$.

Denote by $[\omega]^2$ the set of two-element subsets of $\omega$, and let $\mathcal U$ be an ultrafilter on $[\omega]^2$. Define $\xi\colon \F\To\Pro$ by sending $f\in\F$ to $\pi^2_1$ if and only if the set of all $S\in [\omega]^2$ on which $f$ behaves like the projection to the first coordinate is an element of $\mathcal U$, and to $\pi^2_2$ otherwise. We claim that $\xi$ is a partial clone homomorphism with domain $\F$. To see this, let $f,g,h\in \F$ be given. Let $i,j,k\in\{1,2\}$ be so that $\xi(f)=\pi^2_i$, $\xi(g)=\pi^2_j$, and $\xi(h)=\pi^2_k$. Because $\mathcal U$ is an ultrafilter, and by the definition of $\xi$, the set $Q$ of all $S\in [\omega]^2$ such that $f$ behaves like the $i$-th binary projection, $g$ behaves like the $j$-th binary projection, and $h$ behaves like the $k$-th binary projection on $S$ is an element of $\mathcal U$. Then $f(g(x,y),h(x,y))$ behaves like the composition of those projections on all $S\in Q$, and hence $\xi(f(g(x,y),h(x,y)))=\pi^2_i(\pi^2_j(x,y),\pi^2_k(x,y))=\xi(f)(\xi(g)(x,y),\xi(h)(x,y))$. Therefore, $\xi$ is indeed a partial clone homomorphism, and thus it extends to $\fC$ by Lemma~\ref{lem:extendfrombinary}. 

We claim that $\xi$ is continuous if and only if $\mathcal U$ is principal. If $\mathcal U$ is principal, then there exists $S\in [\omega]^2$ such that  for all $f\in\F$ we have that $\xi(f)=\pi^2_1$ if and only if $f$ behaves like the projection to the first coordinate on $S$; this is a clopen subset of $\F$, and so $\xi$ is continuous. Moreover, we then have that its extension to $\fC$ is continuous by Lemma~\ref{lem:continuousbinary}. Now suppose that $\mathcal U$ is non-principal. Then for any $f\in\F$ and any finite set $A\subseteq \omega^2$, the restriction of $f$ to $A$ can be extended to both a function $f'\in\F$ such that $\xi(f')=\pi^2_1$ and a function $f''\in\F$ such that $\xi(f'')=\pi^2_2$. Hence, $\xi$ is not continuous. 
\end{proof}

We remark that in the preceding proposition, the projective homomorphisms of $\fC$ are precisely those induced by ultrafilters. We can slightly modify this example in order to obtain a \emph{closed} function clone with a discontinuous projective homomorphism, as follows.

\begin{prop}\label{prop:discont2}
There exists a closed function clone of conservative functions on a countable set with a discontinuous projective homomorphism.
\end{prop}
\begin{proof}
Let $\fC$ be the set of all finitary conservative functions on $\omega$ which agree with a projection on every two-element subset of $\omega$. Clearly $\fC$ is a closed function clone. As in the proof of Proposition~\ref{prop:discont1}, any ultrafilter $\mathcal U$ on $[\omega]^2$ defines a projective homomorphism $\xi$: for an $n$-ary function $f\in\fC$, we set $\xi(f):=\pi^n_i$ iff the set of all elements of $[\omega]^2$ on which $f$ behaves like the projection to the $i$-th coordinate is in $\mathcal U$. As before, $\xi$ is continuous if and only if $\mathcal U$ is principal.
\end{proof}

We now present an example of an  \emph{oligomorphic} function clone with a discontinuous projective homomorphism.

\begin{prop}\label{prop:discont3}
There exists a closed oligomorphic function clone on a countable set with a discontinuous projective homomorphism.
\end{prop}
\begin{proof}
We construct the desired clone as the polymorphism clone of a relational structure with a first-order definition in a well-known structure $\Delta$, due to Cherlin and Hrushovski, without the small index property (also see~\cite{Lascar}). 
The signature $\tau$ of $\Delta$ contains
a relation symbol $R_n$ of arity $2n$ for each $n\geq 1$.
The class of all finite $\tau$-structures where each $R_n$ is interpreted as 
an equivalence relation on $n$-tuples of distinct entries with two
equivalence classes is a Fra\"{i}ss\'{e} class.
Let $\Delta$ be its Fra\"{i}ss\'{e} limit, with domain $D$; 
it is $\omega$-categorical
since it is homogeneous and has for all $n\geq 1$ only finitely
many inequivalent atomic formulas with $n$ variables. 

Let $\Gamma$ be the structure with domain $D$ that has for all $n\geq 1$ the relation $R_n$, as well as the $3n$-ary relation
$$S_n:=\big\{ (x,y,z) \in D^{3n} \; \big | \; \neg\big( R_n(x,y)\wedge
 R_n(y,z)\wedge R_n(z,x) \big) \big\} \; .$$
Then $\Gamma$ is first-order definable over $\Delta$ and therefore also
$\omega$-categorical. Since the elements of $\Pol(\Gamma)$ preserve $R_n$ for each $n \geq 1$, the function clone $\Pol(\Gamma)$, viewed as a topological clone, acts naturally on the equivalence classes of $R_n$. Write $\xi_n$ for the mapping which sends every $f\in \Pol(\Gamma)$ to its corresponding function on the equivalence classes of $R_n$. Then $\xi_n$ is a continuous clone homomorphism, and its image is a function clone
on a domain with two elements, which we will denote by $0$ and $1$ in the following (independently of $n$, since the name of the elements of the base set is irrelevant). 

We claim that for every $f \in \Pol(\Gamma)$, the operation $\xi_n(f)$ 
depends on one of its arguments only. 
To see this, observe that
$\xi_n(f)$ preserves the Boolean
relation $\{0,1\}^3 \setminus \{(0,0,0),(1,1,1)\}$
because $f$ preserves $S_n$.
It is well-known that Boolean functions that preserve this Boolean relation depend on one argument only~\cite{Post}. 

Let $\mathcal U$ be a non-principal ultrafilter on $\omega$.
Let $\xi \colon \Pol(\Gamma) \to {\Pro}$ be the mapping which sends every $k$-ary $f\in \Pol(\Gamma)$ to the projection $\pi^k_i\in\Pro$ if and only if the set 
$$\{n\geq 1\;|\; \xi_n(f) \text{ depends on the } i\text{-th argument}\}$$ is an element of $\mathcal U$. Similarly as in the proof of Proposition~\ref{prop:discont1}, one can check that $\xi$ is a clone homomorphism. Moreover,  $\xi$ is not continuous. To see this, observe that for any $S\subseteq \omega$ there exists a binary $f\in\Pol(\Gamma)$ such that $\xi_n(f)$ depends on the first argument if and only if $n\in S$. This function $f$ can be constructed by defining a structure on $D^2$ in the language of $\Delta$ in which for each $2n$-tuple $t\in D^2$ membership in $R_n$ depends only on membership in $R_n$ of the projection of $t$ onto its first coordinate when $n\in S$, and onto its second coordinate when $n\notin S$. Choosing $f$ as any embedding of this structure into $\Delta$ using universality, we obtain a polymorphism of $\Gamma$ with the desired property. But since membership in $\mathcal U$ cannot be determined on any finite subset of $S$, the discontinuity of $\xi$ follows.
\end{proof}

Let us note that the function clone constructed in the preceding proposition also has continuous projective homomorphisms: for example, each single $\xi_n$ is continuous, and the image of $\Pol(\Gamma)$ under $\xi_n$ has a projective homomorphism which is necessarily continuous since the topology on the image is discrete.

\section{Locally Finite Idempotent Algebras}
\label{sect:locfin}

\begin{defn}
We call a function clone $\fC$ \emph{locally finite} iff any algebra which has the functions of $\fC$ as its fundamental operations is locally finite; that is, for all finite subsets $A$ of the domain of $\fC$, the set $\{f(a)\;|\; a \text{ is a tuple of elements in } A \text{ and } f\in \fC\}$ is finite.
\end{defn}

\begin{prop}\label{prop:locfin}
Let $\fC$ be a locally finite idempotent closed function clone with a projective homomorphism. Then $\fC$ has a continuous projective homomorphism.
\end{prop}
\begin{proof}
Let $D$ be the domain of $\fC$. 
Let $A_1\subseteq A_2\subseteq\cdots$ be a sequence of finite subsets of $D$ such that $\bigcup_{n=1}^{\infty} A_n=D$ and such that $A_n$ is closed under the operations in $\fC$ for all $n\geq 1$. For each $n\geq 1$, the restriction of the functions in $\fC$ to $A_n$ induces
a clone $\fC_n$ on $A_n$.

Any projective homomorphism of any $\fC_n$ is clearly continuous, since $\fC_n$ is discrete. 
Moreover, the map from $\fC$ to $\fC_n$ which sends every function in $\fC$ to its restriction to $A_n$ is also continuous and a homomorphism. Therefore, if some $\fC_n$ has a projective homomorphism, then $\fC$ has a continuous projective homomorphism and we are done.

So assume henceforth that no $\fC_n$ has a projective homomorphism. Then, by finite universal algebra, each $\fC_n$ contains a function $f_n$ of arity four satisfying the so-called Siggers identities~\cite{Siggers}. We claim that $\fC$ also has a function which satisfies the Siggers identities. To see this, 
note that for each $n\geq 1$, there is a finite number of possible functions of arity four on $A_n$. Let $T$ be the tree whose vertices on level $n\geq 1$ are precisely the functions of arity four in $\fC_n$ which satisfy the Siggers identities; adjacency between functions of consecutive levels is defined by restriction. 
Then this tree is finitely branching and has vertices on all levels, so by K\H{o}nig's lemma it 
has an infinite branch. The union over the functions of this branch is a function $f$ defined on all of $D$. Clearly, $f$ satisfies the Siggers identities, and $f\in\fC$ since $\fC$ is closed. Since projections do not satisfy the Siggers identities, $\fC$ does not have a projective homomorphism, a contradiction.
\end{proof}

\section{Oligomorphic Clones}
\label{sect:oligo}
\subsection{Canonical function clones}

\begin{defn}
Let $\Delta$ be a structure with domain $D$, and let $f$ be an $n$-ary operation on $D$, where $n\geq 1$. Then $f$ is called \emph{canonical} with respect to $\Delta$ iff for all $k\geq 1$, all tuples $a_1,\ldots,a_n\in D^k$, and all $\alpha_1,\ldots,\alpha_n\in\Aut(\Delta)$ there exists $\beta\in\Aut(\Delta)$ such that
$$
f(\alpha_1(a_1),\ldots,\alpha_n(a_n))=\beta(f(a_1,\ldots,a_n))\;,
$$
where $f$ and $\alpha_i$ are applied componentwise. For $\omega$-categorical $\Delta$, this property has been stated in the language of model theory as `tuples of tuples of the same type are sent to tuples of the same type under $f$'~\cite{RandomMinOps,BP-reductsRamsey, BPT-decidability-of-definability}.
\end{defn}

\begin{defn}
Let $\Delta$ be a structure with domain $D$ and let $k\geq 1$. Denote the set of orbits of the action of $\Aut(\Delta)$ on $k$-tuples by $\mathcal T_k$. Every $n$-ary canonical operation $f$ with respect to $\Delta$ defines an $n$-ary operation $\xi_k^{\typ}(f)$ on $\mathcal T_k$: when $O_1,\ldots,O_n\in\mathcal T_k$, then $\xi_k^{\typ}(f)(O_1,\ldots,O_n)$ is the orbit of $f(o_1,\ldots,o_n)$ where $o_i\in O_i$ can be chosen arbitrarily for $1\leq i\leq n$.
\end{defn}

Consequently, when $\fC$ is a function clone on $D$ consisting of canonical functions with respect to $\Delta$, then $\fC$ defines a set of functions on $\mathcal T_k$, which is easily seen to be a function clone. In fact, the mapping $\xi_k^{\typ}$ which sends every $f\in\fC$ to $\xi_k^{\typ}(f)$ is a continuous clone homomorphism.

\begin{defn}
Let $\Delta$ be a structure with domain $D$ and let $k\geq 1$. For a function clone $\fC$ of canonical functions with respect to $\Delta$, we write $\fC_k^{\typ}$ for the function clone which is the continuous homomorphic image of $\fC$ under $\xi_k^{\typ}$.
\end{defn}

In theory, the clones $\fC_k^{\typ}$ carry more and more information about $\fC$ the larger $k$ gets. However, when $\Delta$ is a homogeneous structure in a finite relational language, then $\fC_k^{\typ}$ and $\fC_m^{\typ}$ are isomorphic for all $k\geq \max\{m,2\}$, where $m$ is the maximal arity of the relations of $\Delta$. 

\begin{defn}
In this situation, we write $\fC_\infty^{\typ}$ for $\fC_m^{\typ}$, and $\xi_\infty^{\typ}$ for $\xi_m^{\typ}$.
\end{defn}

\begin{defn}
Let $\tau$ be a functional signature, and let $\fC$ be a function clone. A set $\Sigma$ of equations over $\tau$ is \emph{satisfiable} in $\fC$ iff there exists a clone homomorphism $\xi$ from the term clone of the completely free $\tau$-algebra into $\fC$ such that $\xi(s)=\xi(t)$ for every equation $(s,t)\in\Sigma$. 

For a set $\F$ of unary functions of $\fC$, we say that $\Sigma$ is satisfiable \emph{modulo $\F$ from the outside} iff there exists a clone homomorphism $\xi$ from the term clone of the completely free $\tau$-algebra into $\fC$ and for every $(s,t)\in\Sigma$ elements $\beta^{s,t}_s,\beta^{s,t}_t\in \F$ such that $\beta^{s,t}_s\circ\xi(s)=\beta^{s,t}_t\circ\xi(t)$.

In both situations, we call $\xi$ a \emph{satisfying clone homomorphism}.
\end{defn}




\begin{prop}\label{prop:equations}
Let $\Delta$ be a homogeneous structure in a finite relational language, and let $\fC$ be a closed function clone of canonical functions with respect to $\Delta$ such that $\fC\supseteq \Aut(\Delta)$. Suppose that a finite set of equations $\Sigma$ is satisfiable in $\fC_\infty^{\typ}$. Then $\Sigma$ is satisfiable in $\fC$ modulo $\overline{\Aut(\Delta)}$ from the outside. Moreover, if $\xi$ is a satisfying clone homomorphism for 
$\fC_\infty^{\typ}$, then the satisfying clone homomorphism for $\fC$ can be chosen to be any $\xi'$ such that $\xi_\infty^{\typ}\circ \xi'=\xi$.
\end{prop}
\begin{proof}
Fix a satisfying clone homomorphism $\xi$ for $\fC_\infty^{\typ}$, and let $\xi'$ be so that $\xi_\infty^{\typ}\circ \xi'=\xi$. Since $\fC_\infty^{\typ}$ is a factor of $\fC$, such a mapping $\xi'$ exists.
Since $\Sigma$ is finite, by adding dummy variables to the terms appearing in $\Sigma$ we may assume that those terms all have the same arity $n\geq 1$.

We first claim that for all finite subsets $A$ of the domain of $\Delta$ there exist $\alpha^{s,t}_s,\alpha^{s,t}_t\in\Aut(\Delta)$ such that $\alpha^{s,t}_s(\xi'(s))$ and $\alpha^{s,t}_t(\xi'(t))$ agree on $A$ for all $(s,t)\in\Sigma$.  To see this, let $u^1,\ldots,u^n\in A^{{|A|}^n}$ be so that for every tuple $v\in A^n$ there exists $1\leq i\leq {|A|}^n$ with $(u^1_i,\ldots,u^n_i)=v$. Let $U^i$ be the orbit of $u^i$ with respect to the componentwise action of $\Aut(\Delta)$. Then $\xi(s)(U^1,\ldots,U^n)=\xi(t)(U^1,\ldots,U^n)$ for all $(s,t)\in\Sigma$. Therefore, $\xi'(s)(u^1,\ldots,u^n)$ and $\xi'(t)(u^1,\ldots,u^n)$ belong to the same orbit, and hence there exist $\alpha^{s,t}_s,\alpha^{s,t}_t\in\Aut(\Delta)$ such that  $\alpha^{s,t}_s(\xi'(s)(u^1,\ldots,u^n))=\alpha^{s,t}_t(\xi'(t)(u^1,\ldots,u^n))$, for all $(s,t)\in\Sigma$. This  proves our claim.

We now provide a standard compactness argument which shows that we can lift the local satisfaction of $\Sigma$ modulo $\overline{\Aut(\Delta)}$ from the outside to the entire domain of $\Delta$; similar arguments are given, for example, in~\cite{BartoPinskerDichotomy, canonical}. 
Let $(A_j)_{j\in\omega}$ be an increasing sequence of finite subsets of the domain of $\Delta$ whose union is the entire domain. For every $j\in\omega$, let 
$$
r_j:=((\alpha^{s,t,j}_s,\alpha^{s,t,j}_t)\;|\; (s,t)\in\Sigma)
$$
be a tuple of length $2\cdot |\Sigma|$ which enumerates the automorphisms whose existence is guaranteed by the above claim for the finite set $A_j$. Now consider the set
$$
\{\gamma\circ r_j\;|\; j\in\omega\text{ and } \gamma\in\Aut(\Delta)\}\;,
$$
where $\gamma$ is applied to the functions in the tuple $r_j$ componentwise. This set is a subset of $\Aut(\Delta)^{2\cdot |\Sigma|}$. It has been shown in~\cite{Topo-Birk} that for all $k\geq 1$, the space $\overline{\Aut(\Delta)}^k$ factored by the equivalence relation where $(\delta_1,\ldots,\delta_k)$ and $(\delta_1',\ldots,\delta_k')$ are identified iff there exists $\gamma\in\Aut(\Delta)$ such that $(\delta_1,\ldots,\delta_k)=(\gamma\circ \delta_1',\ldots,\gamma\circ\delta_k')$ is compact. Hence, the above set has an accumulation point in $\overline{\Aut(\Delta)}^{2\cdot |\Sigma|}$, which we denote by $((\beta^{s,t}_s,\beta^{s,t}_t)\;|\; (s,t)\in\Sigma)$. Clearly, all components of this tuple are elements of $\overline{\Aut(\Delta)}$, and $\beta^{s,t}_s(\xi'(s))=\beta^{s,t}_t(\xi'(t))$ for all $(s,t)\in\Sigma$.
\end{proof}

In the following \red{theorem} we establish a positive answer to Question~\ref{quest:main} for a certain class of function clones which is of great importance in applications -- cf.~ the discussion and references in Section~\ref{sect:results}.

\begin{thm}
Let $\Delta$ be a homogeneous structure in a finite relational language. Let $\fC$ be a closed function clone of canonical functions with respect to $\Delta$ such that $\fC\supseteq \Aut(\Delta)$. If $\fC$ has a projective homomorphism, then so does $\fC_\infty^{\typ}$. In particular, $\fC$ then has a continuous projective homomorphism.
\end{thm}
\begin{proof}
 If  $\fC_\infty^{\typ}$ has no projective homomorphism, then there is a finite set $\Sigma$ of equations (in an arbitrary signature for the functions in $\fC_\infty^{\typ}$) which is not satisfiable in $\Pro$. By Proposition~\ref{prop:equations}, $\fC$ then satisfies $\Sigma$ modulo outside elementary embeddings. Hence, it cannot have a projective homomorphism, which proves the contraposition of the first statement of the theorem. For the final statement, recall that $\fC$ has a continuous homomorphism onto $\fC_\infty^{\typ}$ by our discussion above, and hence composing homomorphisms we obtain that if $\fC_\infty^{\typ}$ has a projective homomorphism, then $\fC$ has a continuous projective homomorphism.
\end{proof}

\subsection{Non-closed function clones}

We now give a negative answer to Question~\ref{quest:main} if we drop the assumption that the function clone be closed.

\begin{prop}
There exists an oligomorphic function clone on a countable domain which has a projective homomorphism, but no continuous one.
\end{prop}
\begin{proof}
Let $(\mathbb Q;<)$ be the rational numbers with the usual order, and let $\fC$ consist of all finitary functions $f$ on $\mathbb Q$ with the following properties:
\begin{itemize}
\item $f\in \Pol(\mathbb Q;<)$;
\item if $f$ is $n$-ary, then there exists an $i\in\{1,\ldots,n\}$, $a\in\mathbb Q$, and $\alpha \in \Aut(\mathbb Q;<)$ such that $f(u)=\alpha(u_i)$ for all $u\in \mathbb Q^n$ with $a< u_j$ for all $1\leq j\leq n$.
\end{itemize}
It is easy to see that $\fC$ is a function clone. Since $\fC$ contains $\Aut(\mathbb Q;<)$, it is oligomorphic.

We can define a homomorphism $\xi\colon \fC\To \Pro$ by sending every $f$ to $\pi^n_i\in {\Pro}$, where $i$ is as above. This homomorphism $\xi$ is not continuous: for every restriction of an $n$-ary function $f$ to a finite set, there exist extensions of this restriction to functions in $\xi^{-1}(\{\pi^n_i\})$, for all $1\leq i\leq n$. 

However, $\xi$ is the unique homomorphism from $\fC$ to ${\Pro}$, since whenever $f$, $i$, and $\alpha$ are as above there exist unary functions $g_1,\ldots,g_n\in\fC$ such that $f(g_1(x_1),\ldots, g_n(x_n))=\alpha(x_i)$ for all $x_1,\ldots,x_n\in\mathbb Q$, and so $f$ has to be sent to $\pi^n_i$ under any homomorphism.
\end{proof}

We remark that the closure of $\fC$ equals $\Pol(\mathbb Q;<)$, and does not possess any projective homomorphism.

\section{Open Problems}
\label{sect:open}
The following questions remain open. 

\begin{quest}
Find a closed function clone where Question~\ref{quest:main} has a negative answer, that is, find
a closed function clone which has a homomorphism to $\Pro$, but no continuous homomorphism to $\Pro$.
\end{quest}

We have mentioned in the introduction that for oligomorphic function clones 
Question~\ref{quest:main} can be reformulated
as a question about the difference between varieties and pseudovarieties. If the
function clone is not oligomorphic, it is not clear whether
the reformulation is still equivalent to Question~\ref{quest:main}.
But the reformulation is of independent interest, in particular in universal algebra, so we explicitly state it here.

\begin{quest}
Let $\mathfrak A$ \red{be an algebra whose operations constitute 
a closed clone over a countably infinite base set. 
Is it true that if} the variety 
generated by $\mathfrak A$ contains a two-element trivial algebra, then so does the pseudovariety generated by 
$\mathfrak A$?
\end{quest}

A positive answer to this question would
imply a positive answer to Question~\ref{quest:main} 
\red{(and the converse \red{is true for oligomorphic clones}; cf.\ Proposition~5 in~\cite{Topo-Birk})}. 

We have seen examples of discontinuous homomorphisms from closed function clones to $\Pro$, but these examples relied on the existence of non-principal ultrafilters. 

\begin{quest}
Is there a model of ZF where \emph{every} homomorphism from a closed function clone to $\Pro$ is continuous?
\end{quest}

Our example of an oligomorphic closed function clone
with a discontinuous homomorphism to $\Pro$ makes essential use of an infinite relational signature.
However, in the context of the constraint satisfaction problem we are particularly interested in finite signatures. 
magenta
\begin{quest}
Let $\Gamma$ be a homogeneous structure with finite relational signature. Is every homomorphism from $\Pol(\Gamma)$ to $\Pro$ continuous?
\end{quest}

\subsection*{\red{Recent Progress}}
After the present article was submitted,  
Barto and Pinsker~\cite{Topo} solved {Question 7.2}  
in the negative; their counterexample is not oligomorphic. Moreover, the importance of uniform continuity (with respect to the metric mentioned in the introduction), rather than continuity, was realized~\cite{uniformbirkhoff, wonderland}; for projective homomorphisms of closed oligomorphic clones, however, this makes no difference. The counterexample mentioned above also provides a counterexample to {Questions 7.1 and 7.3} when continuity is replaced by uniform continuity.

It has been shown that continuity can be dropped in Conjecture~1.2~\cite{Topo}, thus undermining the strategy proposed in the present paper. We believe, however, that Question~1.3 is still of independent mathematical interest, and that its solution could provide valuable insights in connection with  Conjecture~1.2.

An important host of open problems comes from asking analogous questions for minor-preserving maps (also called h1 clone homomorphisms) to $\Pro$, rather than clone homomorphisms to $\Pro$. The significance of minor-preserving maps, in particular for Conjecture~\ref{conj:tractability}, has been recognised in~\cite{wonderland}, and new results in this context can be found in~\cite{TopologyIsRelevant}.

\bibliographystyle{alpha}
\bibliography{homos21.bib}

\newcommand{\etalchar}[1]{$^{#1}$}
\begin{thebibliography}{BMO{\etalchar{+}}19}

\bibitem[BHM10]{BodHilsMartin}
Manuel Bodirsky, Martin Hils, and Barnaby Martin.
\newblock On the scope of the universal-algebraic approach to constraint
  satisfaction.
\newblock In {\em Proceedings of the Annual Symposium on Logic in Computer
  Science (LICS)}, pages 90--99. IEEE Computer Society, July 2010.

\bibitem[Bir35]{Bir-On-the-structure}
Garrett Birkhoff.
\newblock On the structure of abstract algebras.
\newblock {\em Mathematical Proceedings of the Cambridge Philosophical
  Society}, 31(4):433--454, 1935.

\bibitem[BKO{\etalchar{+}}]{BKOPP}
Libor Barto, Michael Kompatscher, Miroslav Ol\v{s}\'{a}k, Trung~Van Pham, and
  Michael Pinsker.
\newblock Equations in oligomorphic clones and the constraint satisfaction
  problem for omega-categorical structures.
\newblock {\em Journal of Mathematical Logic}.
\newblock to appear. Preprint arXiv:1612.07551.

\bibitem[BKO{\etalchar{+}}17]{BKOPP-conf}
Libor Barto, Michael Kompatscher, Miroslav Ol\v{s}\'{a}k, Trung~Van Pham, and
  Michael Pinsker.
\newblock The equivalence of two dichotomy conjectures for infinite domain
  constraint satisfaction problems.
\newblock In {\em Proceedings of the 32nd Annual {ACM/IEEE} Symposium on Logic
  in Computer Science -- LICS'17}, 2017.

\bibitem[BM18]{BodMot-Unary}
Manuel Bodirsky and Antoine Mottet.
\newblock A dichotomy for first-order reducts of unary structures.
\newblock {\em Logical Methods in Computer Science}, 14(2), 2018.

\bibitem[BMM18]{MMSNP}
Manuel Bodirsky, Florent Madelaine, and Antoine Mottet.
\newblock {A universal-algebraic proof of the complexity dichotomy for Monotone
  Monadic SNP}.
\newblock In {\em Proceedings of the Symposium on Logic in Computer Science --
  LICS'18}, 2018.
\newblock Preprint available under ArXiv:1802.03255.

\bibitem[BMO{\etalchar{+}}19]{TopologyIsRelevant}
Manuel Bodirsky, Antoine Mottet, Miroslav Ol\v{s}\'ak, Jakub Opr\v{s}al,
  Michael Pinsker, and Ross Willard.
\newblock Topology is relevant (in the infinite-domain dichotomy conjecture for
  constraint satisfaction problems).
\newblock In {\em Proceedings of the Symposium on Logic in Computer Science --
  LICS'19}, 2019.
\newblock Preprint arXiv:1901.04237.

\bibitem[BMPP]{BMPP16}
Manuel Bodirsky, Barnaby Martin, Michael Pinsker, and Andr{\'{a}}s
  Pongr{\'{a}}cz.
\newblock Constraint satisfaction problems for reducts of homogeneous graphs.
\newblock {\em SIAM Journal on Computing}.
\newblock To appear. Preprint arXiv:1602.05819. A conference version appeared
  in the Proceedings of the 43rd International Colloquium on Automata,
  Languages, and Programming, {ICALP} 2016, pages 119:1--119:14.

\bibitem[Bod07]{Cores-journal}
Manuel Bodirsky.
\newblock Cores of countably categorical structures.
\newblock {\em Logical Methods in Computer Science}, 3(1):1--16, 2007.

\bibitem[BOP18]{wonderland}
Libor Barto, Jakub Opr\v{s}al, and Michael Pinsker.
\newblock The wonderland of reflections.
\newblock {\em Israel Journal of Mathematics}, 223(1):363--398, 2018.

\bibitem[BP11]{BP-reductsRamsey}
Manuel Bodirsky and Michael Pinsker.
\newblock Reducts of {R}amsey structures.
\newblock {\em AMS Contemporary Mathematics, vol. 558 (Model Theoretic Methods
  in Finite Combinatorics)}, pages 489--519, 2011.

\bibitem[BP14]{RandomMinOps}
Manuel Bodirsky and Michael Pinsker.
\newblock Minimal functions on the random graph.
\newblock {\em Israel Journal of Mathematics}, 200(1):251--296, 2014.

\bibitem[BP15a]{BodPin-Schaefer}
Manuel Bodirsky and Michael Pinsker.
\newblock Schaefer's theorem for graphs.
\newblock {\em Journal of the ACM}, 62(3):52 pages (article number 19), 2015.
\newblock A conference version appeared in the Proceedings of STOC 2011, pages
  655--664.

\bibitem[BP15b]{Topo-Birk}
Manuel Bodirsky and Michael Pinsker.
\newblock Topological {B}irkhoff.
\newblock {\em Transactions of the American Mathematical Society},
  367:2527--2549, 2015.

\bibitem[BP16a]{BartoPinskerDichotomy}
Libor Barto and Michael Pinsker.
\newblock The algebraic dichotomy conjecture for infinite domain constraint
  satisfaction problems.
\newblock In {\em Proceedings of the 31th Annual IEEE Symposium on Logic in
  Computer Science -- LICS'16}, pages 615--622, 2016.
\newblock Preprint arXiv:1602.04353.

\bibitem[BP16b]{canonical}
Manuel Bodirsky and Michael Pinsker.
\newblock Canonical functions: a proof via topological dynamics.
\newblock Preprint arXiv:1610.09660, 2016.

\bibitem[BP18]{Topo}
Libor Barto and Michael Pinsker.
\newblock Topology is irrelevant.
\newblock Preprint arXiv:1602.04353. A conference version appeared in the
  proceedings of LICS'16 under the title `The algebraic dichotomy conjecture
  for infinite domain constraint satisfaction problems', 2018.

\bibitem[BPP17]{Reconstruction}
Manuel Bodirsky, Michael Pinsker, and Andr\'{a}s Pongr\'acz.
\newblock Reconstructing the topology of clones.
\newblock {\em Transactions of the American Mathematical Society},
  369:3707--3740, 2017.

\bibitem[BPT13]{BPT-decidability-of-definability}
Manuel Bodirsky, Michael Pinsker, and Todor Tsankov.
\newblock Decidability of definability.
\newblock {\em Journal of Symbolic Logic}, 78(4):1036--1054, 2013.
\newblock A conference version appeared in the Proceedings of LICS 2011, pages
  321--328.

\bibitem[Bul17]{BulatovFVConjecture}
Andrei~A. Bulatov.
\newblock A dichotomy theorem for nonuniform {CSP}s.
\newblock In {\em 58th {IEEE} Annual Symposium on Foundations of Computer
  Science, {FOCS} 2017, Berkeley, CA, USA, October 15-17, 2017}, pages
  319--330, 2017.

\bibitem[GP18]{uniformbirkhoff}
Mai Gehrke and Michael Pinsker.
\newblock Uniform {B}irkhoff.
\newblock {\em Journal of Pure and Applied Algebra}, 222(5):1242--1250, 2018.

\bibitem[Hod97]{Hodges}
Wilfrid Hodges.
\newblock {\em A shorter model theory}.
\newblock Cambridge University Press, Cambridge, 1997.

\bibitem[Las91]{Lascar}
Daniel Lascar.
\newblock Autour de la propri\'et\'e du petit indice.
\newblock {\em Proceedings of the London Mathematical Society}, 62(1):25--53,
  1991.

\bibitem[Pos41]{Post}
Emil~L. Post.
\newblock The two-valued iterative systems of mathematical logic.
\newblock {\em Annals of Mathematics Studies}, 5, 1941.

\bibitem[Sig10]{Siggers}
Mark~H. Siggers.
\newblock A strong {M}al'cev condition for varieties omitting the unary type.
\newblock {\em Algebra Universalis}, 64(1):15--20, 2010.

\bibitem[Zhu17]{ZhukFVConjecture}
Dmitriy Zhuk.
\newblock A proof of {CSP} dichotomy conjecture.
\newblock In {\em 58th {IEEE} Annual Symposium on Foundations of Computer
  Science, {FOCS} 2017, Berkeley, CA, USA, October 15-17, 2017}, pages
  331--342, 2017.

\end{thebibliography}

\end{document}